\documentclass[12pt]{article}

\usepackage{amsmath,amssymb,amscd}
\usepackage{cite}
\usepackage{graphicx,epstopdf}
\usepackage{float}
\usepackage{authblk} 
\numberwithin{equation}{section}
\usepackage[usenames]{color}
\usepackage{color}
\usepackage{colortbl}
\usepackage[linktocpage=true,plainpages=false,pdfpagelabels=false]{hyperref}
\definecolor{linkcolor}{rgb}{0.9,0,0}
\definecolor{citecolor}{rgb}{0,0.6,0}
\definecolor{urlcolor}{rgb}{0,0,1}
\hypersetup{
	colorlinks, linkcolor={linkcolor},
	citecolor={citecolor}, urlcolor={urlcolor}
}

\usepackage{soul}

\def\l{\lambda}
\def\m{\mu}

\def\e{\epsilon}

\def\bp{\begin{proposition}}
\def\ep{\end{proposition}}
\def\bt{\begin{theo}}
\def\et{\end{theo}}
\def\be{\begin{equation}}
\def\ee{\end{equation}}
\def\bl{\begin{lemma}}
\def\el{\end{lemma}}
\def\bc{\begin{corollary}}
\def\ec{\end{corollary}}
\def\pr{\noindent{\bf Proof: }}

\def\bd{\begin{definition}}
\def\ed{\end{definition}}

\def\min{{\rm min\,}}
\def\max{{\rm max\,}}

\def\min{{\rm min\,}}
\def\max{{\rm max\,}}

\newtheorem{theo}{Theorem}[section]
\newtheorem{lemma}{Lemma}[section]
\newtheorem{definition}{Definition}[section]
\newtheorem{corollary}{Corollary}[section]
\newtheorem{proposition}{Proposition}[section]

\begin{document}

\title{Geometry and Singularities of Prony varieties}

\author[1]{Gil Goldman}
\author[2]{Yehonatan Salman}
\author[3]{Yosef Yomdin}

\affil[1,2,3]{Department of Mathematics, The Weizmann Institute of Science, Rehovot 76100, Israel}

\medskip

\affil[1]{Email: gilgoldm@gmail.com}
\affil[2]{Email: salman.yehonatan@gmail.com}
\affil[3]{Email: yosef.yomdin@weizmann.ac.il}

\date{}
\maketitle

\begin{flushright}
    {\it To the Memory of John Mather.}
\end{flushright}

\begin{abstract}

We start a systematic study of the topology, geometry and singularities of the Prony varieties $S_q(\mu)$,
defined by the first $q+1$ equations of the classical Prony system $\sum_{j=1}^d a_j x_j^k = \mu_k, \ k= 0,1,\ldots \ .$

Prony varieties, being a generalization of the Vandermonde varieties, introduced in \cite{arnol1986hyperbolic,kostov1989geometric},
present a significant independent mathematical interest (compare \cite{arnol1986hyperbolic,goldman2018accuracy,kostov1989geometric}).
The importance of Prony varieties in the study of the error amplification patterns in solving
Prony system was shown in \cite{akinshin2015accuracy,akinshin2015Rus,akinshin2017accuracy,akinshin2017error,goldman2018accuracy}.
In \cite{goldman2018accuracy} a survey of these results was given, from the point of view of Singularity Theory.

\smallskip

In the present paper we show that for $q\ge d$ the variety $S_q(\mu)$ is diffeomerphic to an intersection of a certain affine subspace 
in the space ${\cal V}_d$ of polynomials of degree $d$, with the hyperbolic set $H_d$.

On the Prony curves $S_{2d-2}$ we study the behavior of the amplitudes $a_j$ as the nodes $x_j$ collide, and the nodes escape to infinity.

We discuss the behavior of the Prony varieties as the right hand side $\mu$ varies, and possible connections of this problem with J. Mather's 
result in \cite{mather1973solutions} on smoothness of solutions in families of linear systems.

\end{abstract}


\section{Introduction}\label{Sec:Intro}
\setcounter{equation}{0}

This paper is devoted to a detailed study of ``Prony varieties'',
which play important role in some problems of Signal Processing (in particular,
in Fourier reconstruction of ``spike-train signals'' - see Section \ref{Sec:Prony.System} below).
We believe that Prony varieties present a significant independent mathematical interest,
especially from the point of view of Singularity Theory (compare \cite{goldman2018accuracy}).

\smallskip

In particular, in the coarse of our study we provide proofs of most of the results announces in
\cite{goldman2018accuracy}. (However, we keep the present paper independent from \cite{goldman2018accuracy},
and give all the necessary definitions).

\smallskip



\subsection{Prony system}\label{Sec:Prony.System}
\setcounter{equation}{0}

We consider the classical {\it Prony system} of algebraic equations, of the form

\be\label{eq:Prony.system1}
\sum_{j=1}^d a_j x_j^k = \mu_k, \ k= 0,1,\ldots,2d-1.
\ee
Here $d$ and the right hand side $\mu=(\mu_0,\ldots,\mu_{2d-1})$ are assumed to be known,
while $a=(a_1,\ldots,a_d) \in {\mathbb R}^d$ and $x=(x_1,\ldots,x_d) \in {\mathbb R}^d$ are the unknowns to be found.

\smallskip

Prony system appears, in particular, in the problem of moment reconstruction of spike-trains,
that is, of one-dimensional signals $F$ which are linear combinations of $d$ shifted $\delta$-functions:
\be \label{eq:equation.model.delta}
F(x)=\sum_{j=1}^{d}a_{j}\delta\left(x-x_{j}\right).
\ee
We assume that the form (\ref{eq:equation.model.delta}) of signals $F$ is a priori known,
but the specific parameters - the amplitudes $a_j$ and the nodes $x_j$ are unknown.
Our goal is to reconstruct them from $2d$ moments $m_k(F)=\int_{-\infty}^\infty x^k F(x)dx, \ k=0,\ldots,2d-1$,
which are known with a possible error bounded by $\e>0$.

\smallskip

An immediate computation shows that the moments $m_k(F)$ are expressed as $m_k(F)=\sum_{j=1}^d a_j x_j^k$.
Hence our reconstruction problem is equivalent to solving the Prony system (\ref{eq:Prony.system1}), with
$\mu_k=m_k(F).$

\smallskip

We will identify the signal $F$ with the tuple $(a,x)$ of the amplitudes and the nodes of
$F$.
We will assume that the nodes $x$ are pairwise different and ordered: $x_1<x_2<\ldots<x_d$,
and denote the space of the nodes by
	$$ {\cal P}^x_d \cong \Delta_d =  \{(x_1,\ldots,x_d) \in {\mathbb R}^d,\; x_1<x_2<\ldots <x_d\}. $$
Denote the space of the amplitudes by $ {\cal P}^a_d $.
Finally denote by ${\cal P}_d={\cal P}^a_d\times {\cal P}^x_d,$ the parameter space of signals $F$ with $d$ nodes.

The space of the moments $m_0,\ldots,m_{2d-1}$ (or of the right-hand sides $\mu=(\m_0,\ldots,\mu_{2d-1})$
of the Prony system (\ref{eq:Prony.system1})) will be denoted by ${\cal M}_d \cong {\mathbb R}^{2d}$.

\medskip

Prony system appears in many theoretical and applied mathematical problems.
There exists a vast literature on Prony and similar systems - see,
as a very small sample, \cite{prony1795essai,auton1981investigation,batenkov2017accurate,batenkov2016stability,batenkov2013accuracy,beylkin2010approximation,peter2013generalized,peter2011nonlinear,plonka2013many,potts2015Fast}
and references therein.

Some applications of Prony system are of major practical importance, and,
in case when some of the nodes $x_j$ nearly collide, it is well known to present major mathematical difficulties,
in particular, in the context of ``super-resolution problem''
(see \cite{akinshin2015accuracy,akinshin2015Rus,akinshin2017accuracy,akinshin2017error,azais2015spike,batenkov2017accurate,candes2013super,candes2014towards,demanet2015recoverability,demanet2013super,donoho1992superresolution,fernandez2016super,morgenshtern2016super,peter2013generalized}
 as a small sample).

\smallskip

Recent papers \cite{akinshin2015accuracy,akinshin2015Rus,akinshin2017accuracy,akinshin2017error}
deal with the problem of ``error amplification'' in solving a Prony system in the case that the nodes $x_1,\ldots,x_d$
nearly collide. Our approach there is independent of a specific
method of inversion and deals with a possible amplification of the measurements errors,
in the reconstruction process, caused by the geometric nature of the Prony system.

\smallskip

The main observations and results in \cite{akinshin2015accuracy,akinshin2015Rus,akinshin2017accuracy,akinshin2017error}
can be shortly summarized as follows:
{\it the incorrect reconstructions, caused by the measurements noise, are spread along certain algebraic
subvarieties $S_q$ in the parameter space ${\cal P}_d$, which we call the ``Prony varieties'' (see the next section)}.

This important fact allows us to better understand the
geometry of error amplification in solving Prony system: on one hand,
we produce on this basis rather accurate upper and lower bounds for the worst
case reconstruction error. On the other hand, we show that in some cases
this geometric information can be used in order to improve the expected reconstruction accuracy.

\smallskip

A survey of the results on the geometry of error amplification,
obtained in \cite{akinshin2015accuracy,akinshin2015Rus,akinshin2017accuracy,akinshin2017error},
is given in \cite{goldman2018accuracy}. This survey stresses the role of
Singularity Theory in study of Prony inversion, and, in particular,
in study of the Prony varieties, and announces some new results in this direction.
However, these results on the geometry and singularities of the Prony
varieties are stated in \cite{goldman2018accuracy} without proofs.

\smallskip

In the present paper we start a systematic study of the topology,
geometry and singularities of the Prony varieties, providing, in particular,
proofs for most of the results announced in \cite{goldman2018accuracy}.

\subsection{Prony varieties}\label{Sec:Prony.varieties}


\bd\label{def:Prony.leaf}
For $\mu\in {\cal M}_d$ and for $q=0,\ldots, 2d-1,$ the Prony variety $S_q=S_q(\mu)$
is an algebraic variety in the parameter space ${\cal P}_d$,
defined by the first $q+1$ equations of the Prony system (\ref{eq:Prony.system1}):
\be\label{eq:Prony.system22}
\sum_{j=1}^d a_j x_j^k = \mu_k, \ k= 0,1,\ldots,q.
\ee
\ed
Thus the variety $S_q(\mu)$ is completely determined by the
first $q+1$ moments $\mu_k, \ k= 0,1,\ldots,q$,
which are preserved along $S_q(\mu)$. Generically, the dimension of the variety $S_q(\mu)$ is $2d-q-1$. The chain
$$
S_0\supset S_1\supset \ldots \supset S_{2d-2}\supset S_{2d-1}
$$
can be explicitly computed (in principle), from the
known measurements $\mu=(\mu_0,\ldots,\mu_{2d-1})\in {\cal M}_d$.
Notice that $S_{2d-1}(\mu)$ coincides with the set of solutions of the ``full'' Prony system (\ref{eq:Prony.system1}).

\smallskip

If in equations (\ref{eq:Prony.system22}) above we fix the amplitudes $a_j$, we obtain the ``Vandermonde varieties'' 
in ${\cal P}^x_d$, as introduced in \cite{arnol1986hyperbolic,kostov1989geometric}. We expect that Prony varieties, 
being, essentially, the ``fiber spaces'', with the Vandermonde ones as the fibers, share important properties of the
last, described in \cite{arnol1986hyperbolic,kostov1989geometric}.

\smallskip

In our approach to solving Prony systems the Prony varieties $S_q$ serve as an approximation to
the set of possible ``noisy solutions'' of (\ref{eq:Prony.system1}) which
appear for a noisy right-hand side $\mu$. The Prony curve $S_{2d-2}$ is especially prominent in the presentation below.

\bigskip

An important fact, found in \cite{akinshin2015accuracy,akinshin2015Rus,akinshin2017accuracy,akinshin2017error}, is that
{\it if the nodes $x_1,\ldots,x_d$ form a cluster of a size $h\ll 1$, while the measurements error is of order $\e$,
then the worst case error in reconstruction of $S_q$ is of order $\e h^{-q}$. Thus, for smaller $q$,
the varieties $S_q$ become bigger, but the accuracy of their reconstruction becomes better.
The same is true for the accuracy with which $S_q$ approximate noisy solutions of (\ref{eq:Prony.system1}).}
That is, the ``true'', as well as the nosy solutions to the Prony system \eqref{eq:Prony.system1} lie inside an
$\e h^{-q}$-neighborhood of the Prony variety $S_q(\mu)$ calculated from noisy moment measurements $\mu$.

\smallskip

In particular, the worst case error in reconstruction of the solution $S_{2d-1}$ of (\ref{eq:Prony.system1})
is $\sim \e h^{-2d+1},$ while the worst case error in reconstruction of the Prony curve $S_{2d-2}$
is of order $\e h^{-2d+2}.$ That is, the reconstruction of the Prony curve $S_{2d-2}$ is $h$
times better than the reconstruction of the solutions themselves.

\smallskip

Consequently, we can split the solution of (\ref{eq:Prony.system1}) into two steps:
first finding, with an improved accuracy, the Prony curve $S_{2d-2}(\mu)$, and then localizing
{\it on this curve} the solution of (\ref{eq:Prony.system1}).
In particular, in the presence of a certain additional a priori information
on the expected solutions of the Prony system (for example, upper and/or lower bounds on the amplitudes),
it was shown in \cite{akinshin2015accuracy,akinshin2015Rus,akinshin2017accuracy,akinshin2017error,goldman2018accuracy} that the Prony curves
can be used in order to significantly improve the overall reconstruction accuracy.

\medskip

We believe that the results of
\cite{akinshin2015accuracy,akinshin2015Rus,akinshin2017accuracy,akinshin2017error,goldman2018algebraic,goldman2018accuracy}, 
as well as connections with Vandermonde varieties,
justify a detailed algebraic-geometric study of the Prony varieties, and of their singularities.
As above, we refer the reader to \cite{goldman2018accuracy} and references therein for
a survey of these results from the point of view of Singularity Theory.


\smallskip

The paper is organized as follows: in Section \ref{Sec:Explicit.Prony.Param}
our main results and their proofs are presented. This includes a
global algebraic-geometric description of the Prony varieties $S_q(\mu)$.
It is convenient to consider separately the cases $q\le d-1$ and $q\ge d$.
The proofs in the first case are given in Section \ref{Sec:q.less.than.d}
and in the second case in Section \ref{Sec:q.ge.than.d}.

In Section \ref{Sec:Sol.Fns.mu} we informally discuss the behavior of the Prony varieties $S_q(\mu)$ as functions of $\mu$.
This leads, via the previous results, to solving parametric linear systems, and to possible connections of this problem with
J. Mather's result in \cite{mather1973solutions} on smoothness of solutions in families of linear systems. 

In Section \ref{Sec:Prony.curves}
we consider the case of Prony curves $S_{2d-2}(\mu)$ and describe the behavior of
the amplitudes at the nodes collision singularities, and the nodes escape to infinity.

In Section \ref{Sec:two.nodes}, as an illustration, a complete description of the
Prony varieties in the case of two nodes is given.


\section{Global description of Prony varieties}\label{Sec:Explicit.Prony.Param}
\setcounter{equation}{0}

In this section we start a global algebraic-geometric and topological investigation of
the Prony varieties. Our main results are as follows:

\bt\label{thm:sing.Sq1}
For each $\mu \in {\cal M}_d$ and for each $q=0,\ldots,2d-1$ the
Prony variety $S_q(\mu)$ is a smooth submanifold of ${\cal P}_d$ (for $q\ge d$,
possibly empty), and its node projection $S^x_q(\mu)$ is a smooth submanifold of ${\cal P}^x_d$.
\et

It is convenient to separate the cases $q\leq d-1$ and $q\ge d$.
In the first case we have the following result:

\bt\label{thm:sing.Sq2}
For each $q\leq d-1$ and for each $\mu\in {\cal M}_d$ the variety $S_q(\mu)$ satisfies the following conditions:

\smallskip

\noindent 1. $S_q(\mu)$ is non-empty, its dimension is equal to $2d-q-1,$
and the equations $(\ref{eq:Prony.system22})$ are regular at each point $P$ of $S_q(\mu)$,
i.e. the rank of their Jacobian at $P$ is $q+1$.

\smallskip

\noindent 2. The nodes $x_1,\ldots,x_d$ and any $d-q-1$ amplitudes (in particular, $a_{q+2},\ldots,a_d$)
form a global regular coordinate system on $S_q(\mu)$.

\smallskip

\noindent 3. The projection $p_x:S_q(\mu)\to {\cal P}^x_d$ is onto.
It defines $S_q(\mu)$ as a trivial affine bundle, naturally isomorphic
to ${\mathbb R}^{d-q-1}\times {\cal P}^x_d$. In particular, $S_q(\mu)$ is topologically trivial.
\et
The case $q\ge d$ is somewhat more complicated. To state the result we need some definitions and notations. Let

\begin{align}\label{eq:Prony.fol.new.3}
	\begin{split}
		\sigma_0(x_1,\ldots,x_d) &= 1\\
		\sigma_1(x_1,\ldots,x_d) &=-(x_1+\ldots+x_d)\\
		\ldots & \\
		\sigma_k(x_1,\ldots,x_d) &= (-1)^k \sum_{1\le i_1 < i_2 < \cdots < i_k\le d}
		x_{i_1}x_{i_2}\cdots x_{i_k} \\
		\ldots & \\
		\sigma_d(x_1,\ldots,x_d) &=(-1)^d x_1x_2\cdot \ldots \cdot x_{d-1}x_d,
	\end{split}
\end{align}
be the elementary symmetric (Vieta) polynomials in $x_1,\ldots,x_d,$.
Thus $\sigma_j$ are the coefficients of the univariate polynomial

$$
Q(z)=\prod_{j=1}^d(z-x_j)=z^d+\sigma_1z^{d-1}+\ldots+\sigma_d = \sum_{i=0}^d \sigma_{d-i}z^i,
$$
whose roots are the nodes $x_1,\ldots,x_d$.

\smallskip

Let ${\cal V}_d\cong {\mathbb R}^d$ be the space of the coefficients $\sigma=(\sigma_1,\ldots,\sigma_d)$
of the polynomials $Q(z)$ (which we identify with the space of the monic polynomials $Q$ themselves).

\smallskip

Consider a subset $H_d\subset {\cal V}_d$, consisting of {\it hyperbolic} polynomials $Q$,
i.e. of those polynomials $Q(z)=z^d+\sigma_1z^{d-1}+\ldots+\sigma_d$ with all the roots
real and distinct. Thus we consider {\it the open} hypebolic set, excluding the boundary.

\smallskip

The set $H_d$ is important in many problems, and it was intensively studied (see, as a small sample, \cite{arnol1986hyperbolic,kostov1989geometric,kostov2011topics} and references therein).
\bd\label{def:root.vieta.map}
The ``root mapping'' $RM_d:H_d\to {\cal P}^x_d$ is defined by
$$
	RM_d(Q)=x=(x_1,\ldots,x_d)\in {\cal P}^{x}_d,
$$
where $x_1<x_2< \ldots < x_d$ are the ordered roots of the hyperbolic polynomial $Q(z)\in H_d$.

\smallskip

The ``Vieta mapping'' $VM_d: {\cal P}^x_d\to H_d$ is defined by
$$
VM_d(x_1,\ldots,x_d)=\left (\sigma_1(x_1,\ldots,x_d),\ldots, \sigma_d(x_1,\ldots,x_d)\right ),
$$
where $\sigma_i=\sigma_i(x_1,\ldots,x_d),\ i=1,\ldots,d,$ are the Vieta symmetric polynomials in $x_1,\ldots,x_d,$.
\ed
Clearly, on $H_d$ the root mapping $RM_d$ is regular, and $RM_d=VM^{-1}_d$.
Therefore both the mappings $RM_d:H_d\to {\cal P}^x_d$ and its inverse $VM_d: {\cal P}^x_d\to H_d$ provide a regular algebraic diffeomorphism between
$H_d$ and ${\cal P}^x_d$.

\smallskip

Let $\mu \in {\cal M}_d$ be given. For $q\ge d$ consider the following system
of $q-d+1$ linear equations for $\sigma_1,\ldots,\sigma_d$:

\begin{equation}\label{eq:Prony.fol.new.4}
\begin{array}{c}
\mu_{d-1}\sigma_1+\mu_{d-2}\sigma_2+\ldots+\mu_0\sigma_d =-\mu_d\\
\mu_{d}\sigma_1+\mu_{d-1}\sigma_2+\ldots+\mu_1\sigma_d =-\mu_{d+1}\\

..........\\
\mu_{q-1}\sigma_1+\mu_{q-2}\sigma_2+\ldots+\mu_{q-d}\sigma_d =-\mu_{q}\\
\end{array}
\end{equation}
Taking into account that $\sigma_0=1,$ this system can be rewritten as
$$
\sum_{i=0}^d\mu_{l-i}\sigma_i=0, \ l=d,\ldots,q.
$$
For a signal $F$ with nodes $x_1,\ldots,x_d$ and moments $\mu=(\mu_0,\ldots,\mu_q)$, system (\ref{eq:Prony.fol.new.4}) forms a part of the standard
(and classical) linear system for the coefficients of the polynomial $Q$ (see, for instance,
\cite{prony1795essai,peter2013generalized,potts2015Fast}). For $q=2d-1$ the complete system is obtained.

\smallskip

Equations (\ref{eq:Prony.fol.new.4}) define an affine subspace $L_q(\mu)\subset {\cal V}_d$, which is generically
of dimension $2d-q-1$ (but, depending on $\mu$, $L_q(\mu)$ may be empty, or of any dimension not smaller than $2d-q-1$).
We denote by $L^h_q(\mu)$ the intersection of $L_q(\mu)$ and the set $H_d$ of hyperbolic polynomials.

\smallskip

Finally, we notice that system (\ref{eq:Prony.fol.new.4}), being a
linear system in variables $\sigma_1,\ldots,\sigma_d$, forms a
nonlinear system in $x_1,\ldots,x_d,$ if we consider $\sigma_j$ as the Vieta elementary
symmetric polynomials in $x_1,\ldots,x_d$.

\smallskip

Now we have all the tools required to describe the Prony varieties $S_q(\mu)$ for $q\ge d$:

\bt\label{thm:sing.Sq3}
For each $q\geq d$ and for any $\mu\in {\cal M}_d$ the variety $S_q(\mu)$ satisfies the following conditions:

\smallskip

\noindent 1. Either $S_q(\mu)$ is empty, or it is smooth and its dimension is
greater than or equal to $2d-q-1$.

\smallskip

\noindent 2. $S^x_q(\mu)$ is defined in ${\cal P}^x_d$ by system (\ref{eq:Prony.fol.new.4}),
considered as a non-linear system in $x_1,\ldots,x_d$. The Vieta mapping $VM$
and its inverse root mapping $RM$ provide a diffeomorphism between $S^x_q(\mu)$ and $L^h_q(\mu)$.

\smallskip

\noindent 3. The projection $p_x:S_q(\mu)\to S^x_q(\mu),$ as well as its inversion,
are one to one, and provide a diffeomorphism between $S_q(\mu)$ and $S^x_q(\mu)$.
\et
We prove Theorems \ref{thm:sing.Sq1} - \ref{thm:sing.Sq3} in the following order:
first Theorem \ref{thm:sing.Sq2}, then Theorem \ref{thm:sing.Sq3}.
Theorem \ref{thm:sing.Sq1} then follows directly.

\subsection{The case $q\le d-1$. Proof of Theorem \ref{thm:sing.Sq2}}\label{Sec:q.less.than.d}

A direct computation shows that the upper left $d\times d$ minor of
the Jacobian matrix $JP_{2d-1}$ of the Prony system of equations
(\ref{eq:Prony.system1}) is the Vandermonde matrix on the nodes $x_1,\ldots,x_d$.
By the construction, for each signal $F$ in ${\cal P}_d$ these nodes are pairwise different,
and hence the first $d$ rows of $JP_{2d-1}(F)$ are linearly independent.
Therefore, for each $q\le d-1$ all the $q+1$ rows of the Jacobian matrix $JP_{q}(F)$
of the partial system (\ref{eq:Prony.system22}) are linearly independent.
This proves, via Implicit function theorem, that the variety $S_q(\mu)$
is smooth, of dimension $2d-q-1$, which is statement 1 of Theorem \ref{thm:sing.Sq2}.

\smallskip

In order to prove statement 2 we rewrite equations (\ref{eq:Prony.system22}) as
\begin{equation}\label{eq:Prony.fol.new.2}
\begin{array}{c}
a_1+a_2+\ldots+a_{q+1}=\mu_0-a_{q+2}-\ldots - a_d\\
a_1x_1+a_2x_2+\ldots+a_{q+1}x_{q+1}=\mu_1-a_{q+2}x_{q+2}-\ldots - a_dx_d\\
a_1x^2_1+a_2x^2_2+\ldots+a_{q+1}x^2_{q+1}=\mu_2-a_{q+2}x^2_{q+2}-\ldots - a_dx^2_d\\
..........\\
a_1x^{q}_1+a_2x^{q}_2+\ldots+a_{q+1}x^{q}_{q+1}=\mu_{q}-a_{q+2}x^{q}_{q+2}-\ldots - a_dx^{q}_d\\
\end{array}
\end{equation}
The left hand side of (\ref{eq:Prony.fol.new.2}) is the Vandermonde
linear system on the pairwise different nodes $x_1,\ldots,x_{q+1}$ with respect to $a_1,\ldots,a_{q+1}$.
Hence we can uniquely express from (\ref{eq:Prony.fol.new.2})
the amplitudes $a_1,\ldots,a_{q+1}$ via the Cramer rule.
The resulting expressions will be linear in $\mu$ and in $a_{q+2},\ldots,a_d,$
with the coefficients - rational functions in the nodes.
Their denominator is the Vandermonde determinant $V_q(x_1,\ldots,x_{q+1})=\prod_{1\leq i < j \leq q+1}(x_j-x_i)$,
which does not vanish at the points $x=(x_1,\ldots,x_{q+1},\ldots,x_d)\in {\cal P}^x_d$.
Thus (\ref{eq:Prony.fol.new.2}) regularly expresses $a_1,\ldots,a_{q+1}$
through $x_1,\ldots,x_d,a_{q+2},\ldots,a_d$.
We conclude that these last parameters can be considered as a global regular coordinate
system on $S_q(\mu)$. This proves statement 2 of Theorem \ref{thm:sing.Sq2}.

\smallskip

\smallskip

In order to prove statement 3 we notice that by (\ref{eq:Prony.fol.new.2}),
for any fixed $x=(x_1,\ldots,x_d)$ the fiber over $x$ of
the projection $p_x:S_q(\mu)\to {\cal P}^x_d$ is an
affine subspace in ${\cal P}_d$ of dimension $d-q-1,$
regularly parametrized by $a_{q+2},\ldots,a_d$.
Therefore (\ref{eq:Prony.fol.new.2}) provides an
isomorphism of the affine fiber bundles $p_x:S_q(\mu)\to {\cal P}^x_d$
and ${\mathbb R}^{d-q-1}\times {\cal P}^x_d$.
This completes the proof of statement 3 and of Theorem \ref{thm:sing.Sq2}. $\square$

\medskip

Let us stress a special case $q=d-1$. In this case the Prony variety $S_{d-1}(\mu)$ has dimension $d$,
and, according to Theorem \ref{thm:sing.Sq2}, the nodes $x_1<x_2<\ldots <x_d$
can be taken as the coordinates on $S_{d-1}(\mu)$. The Cramer rule applied to (\ref{eq:Prony.fol.new.2}) gives

\be\label{eq:Cramer1}
a_j=\frac{1}{V_d(x_1,\ldots,x_d)} \sum_{l=0}^q A^j_l(x_1,\ldots,x_d)\mu_l, \ j=1,\ldots,d,
\ee
with $A^j_l$ the corresponding minors of the Vandermonde matrix of (\ref{eq:Prony.fol.new.2}).
In fact, the coefficients in (\ref{eq:Cramer1}) can be written in a much simpler form.

\smallskip

By a certain misuse of notations, let us denote by $\varrho_k(u_1,\ldots,u_{d-1}), \ k=1,\ldots,d-1,$
the Vieta symmetric polynomials in $d-1$ variables.

\smallskip

For $x=(x_1,\ldots,x_d)$ we put $\pi_i(x)=(x_1,\ldots,x_{i-1},x_{i+1},\ldots,x_d), \ i=1,\ldots,d,$
and consider the ``partial'' symmetric polynomials
$$
\varrho_k (\pi_i(x)), \ k=1,\ldots,d-1, \ i=1,\ldots,d.
$$
Let $L_{i}(x) =\prod_{1\leq l\le d, \ l\neq i}(x_{i} - x_{l})$ be the denominators in the
summands of the Lagrange interpolation polynomial on the nodes $x_1,\ldots,x_d$.

\bp\label{prop:ampl.fprmula.yoni}
The amplitudes $a_i$ of the points on the Prony variety $S_{d-1}(\mu)$ are
given through the nodes $x=(x_1,\ldots,x_d)$ via the expressions
\be\label{eq:ampl.fprmula.yoni}
a_i=\frac{P(\pi_i(x))}{L_{i}(x)}
\ee
where the polynomial $P$ of $d-1$ variables is defined for $u=(u_1,\ldots,u_{d-1})$ by
$$
P(u)=\mu_0\varrho_{d - 1}(u)+\ldots +\mu_{d-2}\varrho_{1}(u)+\mu_{d-1}.
$$
\ep
\pr
Observe that for $q=d-1$ system (\ref{eq:Prony.fol.new.2}) can be rewritten as follows:

\begin{equation*}
	\begin{pmatrix}
		1 & 1 & ... & 1\\
		x_{1} & x_{2} & ... & x_{d}\\
		& & ........& \\
		x_{1}^{d - 1} & x_{2}^{d - 1} & ... & x_{d}^{d - 1}
	\end{pmatrix}
	\begin{pmatrix}
		a_{1}\\
		a_{2}\\
		...\\
		a_{d}		
	\end{pmatrix}
	=
	\begin{pmatrix}
		\mu_{0}\\
		\mu_{1}\\
		...\\
		\mu_{d - 1}	
	\end{pmatrix}.
\end{equation*}
\vskip0.3cm
\hskip-0.6cm
Using the well known formula for the inverse of the
Vandermonde matrix (see, for example \cite{turner1966inverse}) we obtain from the last matrix equation that
\begin{equation}
	\begin{pmatrix}
		a_{1}\\
		a_{2}\\
		...\\
		a_{d}	
	\end{pmatrix}
	=
	\begin{pmatrix}
		\frac{\varrho_{d - 1}(\pi_{1}(x))}{L_{1}(x)} & \frac{\varrho_{d - 2}(\pi_{1}(x))}{L_{1}(x)} & ... & \frac{\varrho_{1}(\pi_{1}(x))}{L_{1}(x)} & \frac{1}{L_{1}(x)}\\
		\frac{\varrho_{d - 1}(\pi_{2}(x))}{L_{2}(x)} & \frac{\varrho_{d - 2}(\pi_{2}(x))}{L_{2}(x)} & ... &
		\frac{\varrho_{1}(\pi_{2}(x))}{L_{2}(x)} & \frac{1}{L_{2}(x)}\\
		..........\\
		\frac{\varrho_{d - 1}(\pi_{d}(x))}{L_{d}(x)} & \frac{\varrho_{d - 2}(\pi_{d}(x))}{L_{d}(x)} & ... &
		\frac{\varrho_{1}(\pi_{d}(x))}{L_{d}(x)} & \frac{1}{L_{d}(x)}\\
	\end{pmatrix}
	\begin{pmatrix}
		\mu_{0}\\
		\mu_{1}\\
		...\\
		\mu_{d - 1}	
	\end{pmatrix}.
\end{equation}

Therefore for $i=1,\ldots,d$ we get $a_i=\frac{P(\pi_i(x))}{L_{i}(x)}$, with
$$
P(u)=\mu_0\varrho_{d - 1}(u)+\ldots +\mu_{d-2}\varrho_{1}(u)+\mu_{d-1}
$$
defined as above. This completes the proof of Proposition \ref{prop:ampl.fprmula.yoni}. $\square$

\medskip

\bc\label{cor:ampl.Prony.q.ge.d}
For each $q\ge d-1$ equations (\ref{eq:ampl.fprmula.yoni}),
expressing the amplitudes $a_1,\ldots,a_d$ through the nodes $x_1,\ldots,x_d$,
remain valid on the Prony varieties $S_q(\mu)$.
\ec
\pr
By definition, for $q\ge d-1$ we have $S_q(\mu)\subset S_{d-1}(\mu).$ $\square$


\subsection{The case $q\ge d$. Proof of Theorem \ref{thm:sing.Sq3}}\label{Sec:q.ge.than.d}

For each $q\ge d$ the dimension $2d-q-1$ of the Prony varieties $S_q(\mu)$ is strictly
smaller than $d$. Consequently, the projections $S^x_q(\mu)$ of $S_q(\mu)$ onto
the nodes subspace ${\cal P}^x_d$ are proper subvarieties in ${\cal P}^x_d$.

\smallskip

By Corollary \ref{cor:ampl.Prony.q.ge.d} we conclude that for each
$x=(x_1,\ldots,x_d)\in S^x_q(\mu)$ the amplitudes $a=(a_1,\ldots,a_d)$ are uniquely defined,
and given by regular expressions (\ref{eq:ampl.fprmula.yoni}). Thus the projection of $S_q(\mu)$
to $S^x_q(\mu)$ is one-to one and regular, as well as its inverse. This proves statement 3 of Theorem \ref{thm:sing.Sq3}.

\smallskip

Next we concentrate on the Prony varieties $S^x_q(\mu)$, and show that they
are defined in ${\cal P}^x_d$ by system (\ref{eq:Prony.fol.new.4}).
We have to eliminate the amplitudes $a_1,\ldots,a_d$ from the equations
(\ref{eq:Prony.system22}). For this purpose we use a modification of the
classical solution method of the Prony system.

\smallskip

First we show that for $q\geq d$ system (\ref{eq:Prony.system22}) implies system (\ref{eq:Prony.fol.new.4}).
Indeed, for each $l=d,\ldots,q$ we obtain, using (\ref{eq:Prony.system22}), that

$$
\sum_{i=0}^{d}\mu_{l-i}\sigma_i = \sum_{i=0}^{d}\sigma_i\sum_{j=1}^d a_jx^{l-i}_j=
\sum_{j=1}^d a_j \sum_{i=0}^{d}\sigma_i x^{l-i}_j=\sum_{j=1}^d a_jx^{l-d}_jQ(x_j)=0,
$$
since each node $x_j$ is a root of $Q(x)$. In other words, for each $(a,x)\in {\cal P}_d$
satisfying system (\ref{eq:Prony.system22}),
the component $x$ satisfies system (\ref{eq:Prony.fol.new.4}). We conclude that
the projection $S^x_q(\mu)$ of $S_q(\mu)$
onto the nodes subspace ${\cal P}^x_d$ is contained in the zero set of system (\ref{eq:Prony.fol.new.4}).

\smallskip

To prove the opposite inclusion, let us assume that $x=(x_1,\ldots,x_d)$
satisfies system (\ref{eq:Prony.fol.new.4}). We uniquely define the amplitudes $a =(a_1,\ldots,a_d)$
from the Vandermonde linear system, formed by the first $d$ equations of system (\ref{eq:Prony.system22}),
according to expressions (\ref{eq:ampl.fprmula.yoni}). Now we form a signal
$$
F(x)=\sum_{j=1}^d a_j\delta(x-x_j)=(a, x)\in {\cal P}_d,
$$
which by construction satisfies the first $d$ equations of system (\ref{eq:Prony.system22}).
It remains to show that the last $q-d+1$ equations of (\ref{eq:Prony.system22}) are satisfied for $F(x)$.

\smallskip

Consider the rational function $R(z)=\sum_{j=1}^d \frac{a_j}{z-x_j}.$ We have $R(z)=\frac{P(z)}{Q(z)}$
for a certain polynomial $P(z)$ of degree $d-1$ and for
$$
Q(z)=\prod_{j=1}^d(z-x_j)=z^d+\sigma_1z^{d-1}+\ldots+\sigma_d,
$$
where $\sigma_i=\sigma_i(x_1,\ldots,x_d),\ i=1,\ldots,d,$ are, as above, the Vieta
elementary symmetric polynomials in $x_1,\ldots,x_d,$.

\smallskip

Developing the elementary fractions in $R(z)$ into geometric progressions, we get

\be\label{eq:Prony.fol.new.5}
R(z)=\sum_{k=0}^\infty \frac{m_k}{z^{k+1}}, \ \ m_k=m_k(F)=\sum_{j=1}^d a_jx_j^k.
\ee
Therefore, the moments $m_k=m_k(F), \ k=0,1,\ldots,$ given by the left hand side
$\sum_{j=1}^d a_jx_j^k$ of system (\ref{eq:Prony.system22}),
are the Taylor coefficients of the rational function $R(z)=\frac{P(z)}{Q(z)},$ with $P(z)$ of degree $d-1$,
and $Q(z)$ of degree $d$.
Starting with $k=d$ these Taylor coefficients $m_k$ of $R$ are known to satisfy the recurrence relation

\be\label{eq:Prony.fol.new.6}
m_k=-\sum_{s=1}^d \sigma_s m_{k-s},
\ee
$\sigma_l$ being the coefficients of the denominator $Q(z)$ of $R(z)$.
Since by the choice of the amplitudes $a_j$ the first $d$ equations of system (\ref{eq:Prony.system22})
are satisfied, we conclude that $m_k=\mu_k, \ k=0,\ldots,d-1.$

\smallskip

Now we use the assumption that system (\ref{eq:Prony.fol.new.4}) is satisfied. Its equations show that
$\mu_k$ satisfy exactly the same recurrence relation till $k=q$. Since the first $d$ terms are the same,
we conclude that in fact $m_k=\mu_k, \ k=0,\ldots,q.$

\smallskip

This means that the entire system (\ref{eq:Prony.system22}) is satisfied. Consequently,
$F\in S_q(\mu),$ and therefore $x=(x_1,\ldots,x_d)\in S^x_q(\mu) \subset {\cal P}^x_d.$
We conclude that the zero set of system (\ref{eq:Prony.fol.new.4})
is contained in $S^x_q(\mu).$ This completes the proof of the fact that $S^x_q(\mu)$
is defined in ${\cal P}^x_d$ by system (\ref{eq:Prony.fol.new.4}).

\smallskip

We conclude that the Vieta mapping $VM$ transforms the points $x=(x_1,\ldots,x_d)$
of $S^x_q(\mu)$ into the hyperbolic polynomials $Q=VM(x)$  belonging to $L^h_q(\mu)\subset {\cal V}_d$.
Conversely, for each $Q\in L^h_q(\mu)$ its image $RM(Q)$ under the root map belongs to $S^x_q(\mu)$.
Therefore the Vieta mapping $VM$ and its inverse root mapping $RM$ provide a diffeomorphism between $S^x_q(\mu)$
and $L^h_q(\mu)$. This completes the proof of statement 2 of Theorem \ref{thm:sing.Sq3}.

\smallskip

In order to prove statement 1 of Theorem \ref{thm:sing.Sq3} we notice that
$S_q(\mu)\subset {\cal P}_d$ is a diffeomorphic image under
$p_x^{-1}\circ RM$ of $L^h_q(\mu) \subset {\cal V}_d$,
the last being a finite union of open domains in an affine
subspace $L_q(\mu)$ of ${\cal V}_d$. Since $L_q(\mu)$ is
defined by system (\ref{eq:Prony.fol.new.4}) of $q-d+1$ linear equations, $L^h_q(\mu)$
is always smooth and either empty or of dimension not smaller than $2d-q-1$.
The same is true for the diffeomorphic images $S^x_q(\mu)$ and $S_q(\mu)$ of $L^h_q(\mu)$.
This completes the proof of statement 1 of Theorem \ref{thm:sing.Sq3},
and of the entire Theorem \ref{thm:sing.Sq3}. $\square$

\subsection{$S_q(\mu)$ and $L_q(\mu)$ as functions of $\mu$}\label{Sec:Sol.Fns.mu}

In Theorem \ref{thm:sing.Sq3} we do not make any assumption on the rank of linear system (\ref{eq:Prony.fol.new.4}).
It is easy to give examples of a right-hand side $\mu=(\mu_0,\ldots,\mu_q)$
of (\ref{eq:Prony.fol.new.4}) for which the solutions of this system form an
empty set, or an affine subspace $L_q(\mu)$ of any dimension not smaller than $2d-q-1$.
Theorem \ref{thm:sing.Sq3} remains true in each of these cases.
Compare a detailed discussion of the situation for two nodes $(d=2)$ in Section \ref{Sec:two.nodes} below.

\medskip

The possible degenerations of system (\ref{eq:Prony.fol.new.4}) are closely
related to the conditions of solvability of Prony system (see, for example,
Theorem 3.6 of \cite{batenkov2013geometry}, and the discussion thereafter).
Both these questions are very important in the robustness analysis of the Prony inversion,
but we do not discuss them here. In Section \ref{Sec:two.nodes}
we illustrate the results above providing a complete description of the Prony curves in the case of two nodes.

\medskip

The observations above lead to a very important question: what can be said about the behavior of the affine subspace $L_q(\mu)$ as a 
function of $\mu$? Via the results above, answering this question will describe also (up to intersection with $H_d$) the behavior 
of the Prony varieties $S_q(\mu), \ q\ge d,$ 
as a function of $\mu$. A very important special case is $q=2d-1$ where $S_{2d-1}(\mu)$ is the set of solutions $F_\mu$ of the original 
Prony system, while $L_{2d-1}(\mu)$ is the set of polynomials $Q_\mu$ whose roots are the nodes of the solutions $F_\mu$.
  
Linear system (\ref{eq:Prony.fol.new.4}) essentially presents a family, parametrized by the moments $\mu=(\mu_0,\mu_1,\ldots)$. 
As it was mentioned above, the rank of this system typically changes with $\mu$, and the behavior of the affine subspaces $L_q(\mu)$, 
as a function of $\mu$, may be rather complicated. 

\smallskip

We expect that J. Mather's theorem in \cite{mather1973solutions}, on smoothness of solutions of parametric families of linear systems, will 
be important in analysis of this problem. Notice, however, that system (\ref{eq:Prony.fol.new.4}) is rigidly structured: its matrices $M(\mu)$ are 
of Hankel type. Consequently, the transversality of the family $\mu\to M(\mu)$ to the rank stratification of the space of Hankel matrices, required
in J. Mather's theorem, is not self-evident. Also the second condition of this theorem, the existence of a solution for each $\mu$, is a delicate 
question, closely related to solvability conditions for the Prony system.

We believe that approaching these problems with the tools used in the proof of J. Mather's theorem, may be very productive.


\subsection{The case $q=2d-2$: Prony curves}\label{Sec:Prony.curves}

The case $q=2d-2$ is especially important in study of error amplification
(see \cite{akinshin2015accuracy,akinshin2015Rus,akinshin2017accuracy,akinshin2017error}).
For generic $\mu \in {\cal M}_d$ the dimension of the Prony variety $S(\mu)=S_{2d-2}(\mu)$ is one,
and by Theorem \ref{thm:sing.Sq3} the variety $S(\mu)$ is a smooth curve consisting of a finite number of open intervals.
These intervals are parametrized, via the root mapping $RM$,
by the intervals of the intersection $L^h_{2d-2}(\mu)$ of the straight line $L_{2d-2}(\mu)$
with the hyperbolic set  $H_d$ in the polynomial space ${\cal V}_d$.

\smallskip

In turn, a convenient explicit parametrization of the straight line $L_{2d-2}(\mu)$
can be obtain as follows: consider system (\ref{eq:Prony.fol.new.4}), with $q=2d-2$,
whose equations define the affine space $L_{2d-2}(\mu)$ in the polynomial space ${\cal V}_d$.
We complete this system to

\be\label{eq:param.S}
\left(\begin{array}{ccccc}
\mu_{0} & \mu_{1} & ... & \mu_{d - 2} & \mu_{d - 1}\\
\mu_{1} & \mu_{2} & ... & \mu_{d - 1} & \mu_{d}\\
& .........\\
\mu_{d - 1} & \mu_{d} & ... & \mu_{2d - 3} & \mu_{2d - 2}\\
\end{array}\right)
\left(\begin{array}{c}
\sigma_{d}\\
\sigma_{d - 1}\\
........\\
\sigma_{1}\\
\end{array}\right)
 = - \left(\begin{array}{c}
\mu_d \\
\mu_{d+1}\\
...\\
\mu_{2d-1}
\end{array}\right).
\ee
which is system (\ref{eq:Prony.fol.new.4}) with $q=2d-1.$
In other words, we add the last equation, corresponding to $q=2d-1.$
The matrix on the left hand side of (\ref{eq:param.S})
is called the moment Hankel matrix $M_d(\mu)$.
This matrix plays the central role in study of Prony systems.

\smallskip

We denote the determinant of $M_d(\mu)$ by $\delta(\mu)$,
and its minimal singular value by $\eta(\mu)$.
Notice that, in fact, the last moment $\mu_{2d-1}$ does not enter $M_d(\mu)$,
and therefore this matrix is constant along the Prony curve $S(\mu)$.
We conclude that for $\delta(\mu)\ne 0$ the rank of the first $d-1$
equations of (\ref{eq:param.S}) is $d-1$, and hence $L_{2d-2}(\mu)$
is one-dimensional, i.e. a straight line in ${\cal V}_d$.

\smallskip

All the entries $\mu_j$ in (\ref{eq:param.S}) besides the bottom
moment $\mu_{2d-1}$ on the right-hand side are fixed on the Prony curve
$S(\mu)$. Accordingly, to obtain a parametrization of $L_{2d-2}(\mu)$
we put $t=\mu_{2d-1}$ and take it as a free parameter.
For each given $\mu_{2d-1}=t$ we solve (\ref{eq:param.S})
and obtain the corresponding coordinates $\sigma_1(t),\ldots, \sigma_d(t)$
of the point $Q_t(z)\in L_{2d-2}(\mu)$. Explicitly, via Cramer's rule and Theorem \ref{thm:sing.Sq3} we have

\bp\label{Prop:param.Prony.curves}
Let $\mu\in {\cal M}_d$ be such that $\delta(\mu)\ne 0$.
Then the line $L_{2d-2}(\mu)\subset {\cal V}_d$ possesses a
parametrization $Q_t(z)=z^d+\sigma_1(t)z^{d-1}+\ldots +\sigma_d(t)$ with $t=\mu_{2d-1}$ and
$$
\sigma_{k}(t) = \alpha_{k}t+\beta_{k}, \ k = 1,...,d, \ \ \text \ {with} \ \ \alpha_k=\frac{(- 1)^l M_{d,l}(\mu)}{\delta(\mu)},
$$
$$
\beta_{k} = \frac{(- 1)^{l}}{\delta(\mu)}\left(\mu_{d}\cdot M_{1, l}(\mu) - \mu_{d + 1}\cdot M_{2, l}(\mu) + ... + (- 1)^{d - 2}\mu_{2d - 2}\cdot M_{d - 1,l}(\mu)\right),
$$
where $l=d-k+1,$ and $M_{i,j}(\mu)$ denotes the minor of the entry
in the $i$-th row and $j$-th column of $M_d(\mu)$.
\ep
Notice that by the construction the roots of $Q_t(z)$ are
always the nodes $x_{1},...,x_{d}$ in the solution $F(t)$ of
the original Prony system (\ref{eq:Prony.system1}),
with the right hand side $(\mu_0,\ldots,\mu_{2d-2},t)$.

\smallskip

To get real and distinct nodes $x_{1},...,x_{d}$ we take only
those values of $t$ for which $Q_t(z)\in H_d$. Thus,
we define $A_{\mu}\subseteq\Bbb R$ as the set of all $t\in\Bbb R$
for which $Q_t(z)\in H_d$. $A_{\mu}$ is a finite union of open intervals in $\Bbb R$.

\smallskip

It is important to accurately describe the behavior of the nodes
and of the amplitudes of $F(t)$ along the Prony curve $S(\mu),$
in terms of the {\it known ``measurements'' $\mu$}. It is explained
in \cite{goldman2018accuracy} how this information can help to improve reconstruction accuracy.
The following two results in this direction were announced in \cite{goldman2018accuracy} without proof:

\smallskip

\noindent 1. Assume that $\delta(\mu)\ne 0$. Then if two nodes $x_i(t),x_{i+1}(t)$
of $F(t)$ collide as $t\to t_0$, then both the amplitudes $a_i(t), a_{i+1}(t)$ tend to infinity as $t\to t_0$.

\smallskip

\noindent 2. Assume that $\delta(\mu)\ne 0$, and that the upper left $(d-1)\times (d-1)$
minor of $M_d(\mu)$ is also non-degenerate. Then for $t\to \pm \infty$ at most one
node of $F(t)$ can tend to infinity.

These results will follow from significantly more accurate results of
Theorem \ref{thm:ampl.Prony.curves}, Theorem \ref{thm:nodes.Prony.curves},
and Proposition \ref{Prop:real.roots} which we prove below.

\subsubsection{Behavior of the amplitudes as the nodes near-collide}\label{Sec:Nodes.Collision.Ampl}

In study of the nodes collisions it is convenient to slightly change the
initial setting of the problem, and to consider {\it unordered, and possibly colliding}
nodes $x=(x_1,\ldots,x_d)\in {\mathbb R}^d$.

\smallskip

Denote by $\Sigma$ the ``diagonal'' in ${\mathbb R}^d$ consisting of
all $x=(x_1,\ldots,x_d)\in {\mathbb R}^d$ with at least two coordinates equal.
For $s=2,\ldots,d$ let $J_s$ be the set of all the $s$-tuples
$$
J=\{1\le r_1< \ldots < r_s \le d\}
$$
We denote by $J^i_s$ the subset of $J_s$ consisting of the $s$-tuples $J$ with $i\in J$.

\smallskip

For $J\in J_s$ let
$$
\Sigma(J)=\{x=(x_1,\ldots,x_d)\in {\mathbb R}^d, \ x_{r_1}=x_{r_2}=\ldots=x_{r_s}\}.
$$
be the subspace in ${\mathbb R}^d$ where the corresponding $s$ coordinates coincide. Then we have

\be\label{eq:diag}
\Sigma=\cup_{s,J\in J_s} \ \Sigma(J).
\ee
We also define the set of all the $s$-collisions as $\Sigma_s=\cup_{J\in J_s} \ \Sigma(J)$,
and the set of all the $s$-collisions, which include the $i$-th coordinate,
as $\Sigma^i_s=\cup_{J\in J^i_s} \ \Sigma(J)$.

\smallskip

Previously we considered the Prony varieties $S^x_q(\mu)$ only inside
the pyramid $\Delta_d$, which is one of the components of ${\mathbb R}^d\setminus
\Sigma$.
But equations (\ref{eq:Prony.fol.new.4}) define $S^x_q(\mu)$ in the entire space ${\mathbb R}^d$,
and by Theorem \ref{thm:sing.Sq2} and Theorem \ref{thm:sing.Sq3} we conclude that in ${\mathbb R}^d$
all the singularities of $S^x_q(\mu)$ are contained in the diagonal $\Sigma$.

\smallskip

Notice that the points of $S^x_q(\mu)\cap \Sigma$ may be non-singular points of $S^x_q(\mu)$:
compare Section \ref{Sec:two.nodes} below. Another remark is that the nodes permutations preserve
the Prony varieties, and hence the study of their ``out-of-collisions''
part can be restricted to the pyramid $\Delta_d$.

\smallskip

Now we consider Prony curves $S(\mu)=S_{2d-2}(\mu)$ and for $F=(a,x)\in S(\mu)$
describe the behavior of the amplitudes $a=(a_1(x),\ldots,a_d(x))$ as
the nodes $x=(x_1,\ldots,x_d)\in S^x(\mu)$ approach the diagonal $\Sigma_d$.

\smallskip

For given $i$ and $s$ with $1\le i,s \le d$ and for
$x=(x_1,\ldots,x_d)\in {\mathbb R}^d$ let $dist^i_s(x)$ denote
the distance of $x$ to the stratum $\Sigma^i_s$ of the diagonal $\Sigma$ (we work with the Euclidean norm in ${\mathbb R}^d$).

Essentially, $dist^i_s(x)$ estimates the minimal size of a cluster of $s$
nodes $x_j$, containing the node $x_i$. Indeed, we have the following result:

\bl\label{lem:cluster.size}
For given $i$ and $s$ with $1\le i,s \le d$ and for
any $x=(x_1,\ldots,x_d)\in {\mathbb R}^d$ there are pairwise
different, and different from $i$, indices $l_1,\ldots,l_{s-1}$ such that
$$
|x_{l_j}-x_i|\le 2 \ dist^i_s(x), \ j=1,\ldots,s-1.
$$
\el
\pr
By definition, there is a point $\bar x=(\bar x_1,\ldots,\bar x_d)\in \Sigma^i_s$
such that $||x-\bar x||=dist^i_s(x)$. In particular, we have $|x_j-\bar x_j|\leq dist^i_s(x), \ j=1,\ldots,d.$
By definition, there are indices
$J=\{1\le r_1< \ldots < r_s \le d\}$,
one of them is equal to $i$, such that $\bar x_{r_1}=\bar x_{r_2}=\ldots=\bar x_{r_s}$.
Re-denote by $l_j, \ j=1,\ldots,s-1,$ the indices $r_m$, different from $i$. Then for each $l_j$ we have
$$
|x_{l_j}-x_i|\le |x_{l_j}-\bar x_{l_j}|+|x_i-\bar x_i|\le 2 \ dist^i_s(x),
$$
since $\bar x_{l_j}=\bar x_i$. This completes the proof. $\square$

\bigskip

Let us recall that by Proposition \ref{prop:ampl.fprmula.yoni} the
amplitudes $a_i$ on $S^x(\mu)$ are uniquely expressed through $x$ as
\be\label{eq:amplit.exp}
a_i(x)=\frac{P(\pi_i(x))}{L_{i}(x)},
\ee
where for $u=(u_1,\ldots,u_{d-1})$ the polynomial $P$ is defined by
\be\label{eq:amplit.exp1}
P(u)=\mu_0\varrho_{d - 1}(u)+\ldots +\mu_{d-2}\varrho_{1}(u)+\mu_{d-1},
\ee
and $L_{i}(x) =\prod_{1\leq l\le d, \ l\neq i}(x_{i} - x_{l})$.

\bigskip

Finally, for $d,\mu$ as above, and for a given $D>1$, put
$$
C_1=\frac{1}{2^dD^{2d-s-1}{\sqrt d}}, \ C_2=\frac{1}{D^{d-1}{\sqrt d}}, \ C_3=\nu_{d-1}(1+D)^{d-1},
$$
with $\nu_{d-1}=\max_{k=0,\ldots,d-1} |\mu_k|$.

\smallskip

Now we have all the definitions and preliminary facts required in
order to state and prove our result on the behavior of the amplitudes
at the nodes collision point on the Prony curve.

\bt\label{thm:ampl.Prony.curves}
Let $\mu\in {\cal M}_d$ be such that $\delta(\mu)\ne 0$.
Then for each $x=(x_1,\ldots,x_d)\in S^x(\mu)$, with $|x_j|\le D, \ j=1,\ldots,d$,
for each $i=1,\ldots,d,$ and for each $s=2,\ldots,d$ we have
\be\label{eq:asympt.ampl}
\frac{C_1\eta(\mu)}{(dist^i_s(x))^s}\le \frac{C_2\eta(\mu)}{L_{i}(x)} \le |a_i(x)| \le \frac{C_3}{L_{i}(x)}.
\ee
In particular, if $x\in S^x(\mu)$ tends to $\bar x \in \Sigma(J)$ with $J\in J_s$,
then for each $i\in J$ the amplitude $a_i(x)$ tends to infinity at least as fast as
$\frac{C_1\eta(\mu)}{||x-\bar x||^s}$.
\et
\pr
Starting with (\ref{eq:amplit.exp}) we reduce bounding $a_i(x)$ to estimating the
numerator and the denominator in this expression. For fixed $i$ and $s$,
and for any $x=(x_1,\ldots,x_d)\in {\mathbb R}^d,\ |x_j|\le D$, (not only for $x\in S^x(\mu)$)
we have the following inequality:
\be\label{eq:asympt.ampl.1}
|L_{i}(x)|\le (2D)^{d-s}[2 \ dist^i_s(x)]^s=2^dD^{d-s}(dist^i_s(x))^s.
\ee
Indeed, by Lemma \ref{lem:cluster.size}, some $s-1$ factors in $L_{i}(x)$ do
not exceed $2 \ dist^i_s(x)$, while the remaining factors are bounded by $2D$.
Accordingly, it remains only to bound the polynomial $P(\pi_i(x))$
from above and from below, for any $x\in S^x(\mu)$.

\smallskip

From (\ref{eq:amplit.exp1}) and from the assumption that $|x_j|\le D$ we get,
denoting, as above, by $\nu_{d-1}$ the maximum $\max_{0\le k \le d-1} \ |\mu_k|$,

\be\label{eq:asympt.ampl.2}
|P(\pi_i(x))|\le \nu_{d-1} \sum_{k=0}^d(^{d-1}_{\ k})D^k=\nu_{d-1}(1+D)^{d-1}.
\ee
To prove the lower bound for $P(\pi_i(x))$ on the Prony curve $S^x(\mu)$,
we add $P$ to the system of equations (\ref{eq:Prony.fol.new.4}),
defining $S(\mu)$, and transform the resulting system to a convenient form.

\smallskip

Let us start with a simple identity for the symmetric polynomials: we fix $i$ and denote $\pi_i(x)$ by $x^*$. Then
$$
\sigma_k(x)=-\varrho_{k-1}(x^*)x_i+\varrho_k(x^*), \ k=1,\ldots,d.
$$
Recall that $\varrho_0\equiv 1, \ \varrho_d\equiv 0.$ We write shorty
$\varrho_k$ for $\varrho_{k}(x^*)$. System (\ref{eq:Prony.fol.new.4}),
defining, for $q=2d-2$, the Prony curve $S^x(\mu)=S^x_{2d-2}(\mu),$ can be now rewritten as $(^*)$:
$$
\begin{array}{c}
\mu_0(\varrho_{d-1}x_i-\varrho_d)+\mu_1(\varrho_{d-2}x_i-\varrho_{d-1})+\ldots+
\mu_{d-1}(\varrho_0x_i-\varrho_1)-\mu_d=0\\
\mu_1(\varrho_{d-1}x_i-\varrho_d)+\mu_2(\varrho_{d-2}x_i-\varrho_{d-1})+
\ldots+\mu_{d}(\varrho_0x_i-\varrho_1)-\mu_{d+1}=0\\

..........\\
\mu_{d-2}(\varrho_{d-1}x_i-\varrho_d)+\mu_{d-1}(\varrho_{d-2}x_i-\varrho_{d-1})+\ldots+\mu_{2d-3}
(\varrho_0x_i-\varrho_1)-\mu_{2d-2}=0\\
\end{array}
$$
Let us assume now that $P(x^*)=\gamma,$ or
$$
\mu_{0}\varrho_{d - 1} + \mu_{1}\varrho_{d - 2} + ... + \mu_{d - 2}\varrho_{1} + \mu_{d - 1}\varrho_0=\gamma
$$
Multiplying this equation by $x_{i}$ and subtracting it from the
first equation in the system $(^*)$, we get all the terms containing the product with $x_{i}$
cancelled, and obtain a new equation
$$
\mu_{1}\varrho_{d - 1} + \mu_{2}\varrho_{d - 2} + ... + \mu_{d - 1}\varrho_{1} + \mu_{d}\varrho_0=\gamma x_{i}.
$$
Multiplying this new equation by $x_{i}$ and subtracting it from the second equation in the system $(^*)$,
we get all the terms containing the product with $x_{i}$ cancelled, and obtain the next equation
$$
\mu_{2}\varrho_{d - 1} + \mu_{3}\varrho_{d - 2} + ... + \mu_{d}\varrho_{1} + \mu_{d+1}\varrho_0=\gamma x^2_{i}.
$$

\medskip

\smallskip

Continuing in this way we obtain the following system:
\vskip-0.2cm

$$\left(\begin{array}{ccccc}
\mu_{0} & \mu_{1} & ... & \mu_{d - 2} & \mu_{d - 1}\\
\mu_{1} & \mu_{2} & ... & \mu_{d - 1} & \mu_{d}\\
& .........\\
\mu_{d - 1} & \mu_{d} & ... & \mu_{2d - 3} & \mu_{2d - 2}\\
\end{array}\right)
\left(\begin{array}{c}
\varrho_{d - 1}\\
\varrho_{d - 2}\\
........\\
\varrho_{1}\\
1
\end{array}\right)
 = \left(\begin{array}{c}
\gamma \\
\gamma x_{i}\\
...\\
\gamma x^{d-2}_{i}\\
\gamma x^{d-1}_{i}
\end{array}\right).$$
Since, by our assumption, the minimal singular value of the matrix $M_d(\mu)$
in the left hand side of the last system is $\eta(\mu)$, we conclude that the
norm of the right hand side is at least $\eta(\mu)$:
\be\label{eq:asympt.ampl.3}
|\gamma|(\sum_{j=0}^{d-1} x^{2j}_{i})^{\frac{1}{2}}\ge \eta(\mu), \ \text {or} \
|\gamma|\ge \frac{\eta(\mu)}{(\sum_{j=0}^{d-1} x^{2j}_{i})^{\frac{1}{2}}}\ge \frac{\eta(\mu)}{D^{d-1}\sqrt d}.
\ee
The last inequality follows since by assumptions $|x_{i}|\le D$. Thus by inequalities
(\ref{eq:asympt.ampl.1}), (\ref{eq:asympt.ampl.2}) and (\ref{eq:asympt.ampl.3}) we have
$$
|L_{i}(x)|\le 2^dD^{d-s}(dist^i_s(x))^s, \ \ \ \frac{\eta(\mu)}{D^{d-1}\sqrt d}\le
|P(\pi_i(x))|\le \nu_{d-1}(1+D)^{d-1},
$$
which implies, via (\ref{eq:amplit.exp})
$$
\frac{\eta(\mu)}{2^dD^{2d-s-1}{\sqrt d} (dist^i_s(x))^s}\le
\frac{\eta(\mu)}{D^{d-1}{\sqrt d} L_{i}(x)} \le |a_i(x)| \le \frac{\nu_{d-1}(1+D)^{d-1}}{L_{i}(x)}.
$$
This completes the proof of Theorem \ref{thm:ampl.Prony.curves}.  $\square$

\subsubsection{Escape of the nodes to infinity}\label{Sec:Nodes.Escape}

In Proposition \ref{Prop:param.Prony.curves} above we describe a natural
parametrization of the Prony curves $S^x(\mu)$ with $M_d(\mu)$ non-degenerate.
It goes via the root mapping $RM$ on the intersection $L^h(\mu)$ of the line $L(\mu)$
with the hyperbolic set $H_d\subset {\cal V}_d$. In turn, the polynomials
$Q_t(z)=z^d+\sigma_1(t)z^{d-1}+\ldots +\sigma_d(t)$ in the line $L(\mu)$ are given by
\be\label{eq:param}
\sigma_{k}(t) = \alpha_{k}t+\beta_{k}, \ k = 1,...,d,
\ee
with $t=\mu_{2d-1}$. The explicit expressions through $\mu_0,\ldots,\mu_{2d-2}$
for $\alpha_k,\beta_k$ are given in Proposition \ref{Prop:param.Prony.curves}.

\smallskip

The set $A_{\mu}\subseteq\Bbb R$ was defined as the set of all $t\in\Bbb R$
for which $Q_t(z)\in H_d$. $A_{\mu}$ is a finite union of open intervals in $\Bbb R$.
The expression $(x_1(t),\ldots,x_d(t))=RM(Q_t)$ provides a diffeomorphic parametrization of the
Prony curve $S^x(\mu)$, for $t\in A_{\mu}$. See Section \ref{Sec:two.nodes}
for examples in the case of two nodes.

\bt\label{thm:nodes.Prony.curves}
Let $\mu\in {\cal M}_d$ be such that the moment Hankel matrix $M_d(\mu)$
is non-degenerate, as well as its top-left $(d-1)\times (d-1)$ minor.
Then for $t\to \pm \infty$ inside $A_\mu$ at most one node among $x_1(t),\ldots,x_d(t)$ can tend to infinity.
\et
\pr
By our assumptions, the parametrization (\ref{eq:param}) above is applicable, and,
by Proposition \ref{Prop:param.Prony.curves}, we have $\alpha_1\ne 0$ in (\ref{eq:param}).
Hence the required result is implied directly by the following statement
(where we do not insist on all the roots of $Q_t(z)$ being real):

\bp\label{Prop:real.roots}
Let $Q_t(z)=z^d+\sigma_1(t)z^{d-1}+\ldots +\sigma_d(t), \ \ \sigma_{k}(t) = \alpha_{k}t+\beta_{k},$
be a polynomial pencil, with $a_1\ne 0$. There are positive constants $t_0,t_1,\l_0,A_1<A_2$
(defined in the proof below) such that for $\alpha_1>0$ and for $t\ge t_0$ the polynomial $Q_t(z)$
has no real roots on the interval $[\l_0,\infty)$.

For $\alpha_1<0$ and for $t\ge t_0$ the polynomial $Q_t(z)$ has no real roots on the
intervals $[\l_0,A_1t]$ and $[A_2t,\infty)$ and exactly one real root on the interval $[A_1t,A_2t]$.
\ep
\pr
We use a version of the Descartes rule, usually called the Budan-Fourier Theorem.
Let $\nu_{Q}(\lambda)$ denote the number of sign changes in the sequence
$(Q(\lambda),Q'(\lambda),...,Q^{(n)}(\lambda))$
(the number of sign changes in a sequence of real numbers is counted with zeroes omitted).

\bl\label{lem:bud.Four}
For any polynomial $Q(z)$ of degree $d$ and $a<b$ the number of
zeros of $Q$ (counted with their multiplicities) in the interval $(a,b]$
is less than or equal to $\nu_{Q}(a) - \nu_{Q}(b)$, and differs from it by an even number.
\el
For the derivatives $Q^{(r)}_t(\l), \ r=0,1,\ldots,d,$ of $Q_t(\l)$
we have, denoting for $n\ge m$ by $A(n,m)$ the product $(n-m+1)(n-m+2)\ldots (n-1)n$,
$$
Q^{(r)}_t(\l)=A(d,r)\l^{d-r}+\sum^{d-r}_{k=1}A(d-k,r)(\alpha_kt+\beta_k)\l^{d-k-r}.
$$
To make the computations more transparent, we normalize the derivatives,
and estimate the signs of $\hat Q^{(r)}_t(\l)=\frac{Q^{(r)}_t(\l)}{A(d,r)\l^{d-r}}$, for which we obtain

\be\label{eq:der.expression}
\hat Q^{(r)}_t(\l)=1+\frac{C(r)t\alpha_1}{\l}\left (1+\frac{\beta_1}{\alpha_1t}+
\sum^{d-r}_{k=2}\frac{B(k,r)}{\l^{k-1}}\left (\frac{\alpha_k}{\alpha_1}+\frac{\beta_k}{\alpha_1t}\right )\right ),
\ee
where $C(r)=\frac{A(d-1,r)}{A(d,r)}=\frac{d-r}{d},$ \ $B(k,r)=\frac{A(d-k,r)}{A(d-1,r)},$
for $r=0,\ldots,d-1$. Notice that $C(d)=0$, and hence $\hat Q^{(d)}_t(\l)=1$. Recall also that by assumptions $\alpha_1\ne 0$.

\bigskip

Next we assume that
$$
t\ge t_0:=20 \ \max_{k=1,\ldots,d}|\frac {\beta_k}{\alpha_1}|,\ \ \l \ge \l_0:=20d \
\max_{k=2,\ldots,d, \ r=0,\ldots,d} B(k,r)(|\frac{\alpha_k}{\alpha_1}|+\frac{1}{20}).
$$

\bl\label{lem:simple.expr}
Under these assumptions we have
\be\label{eq:der.expression1}
\hat Q^{(r)}_t(\l)=1+\frac{C(r)t\alpha_1}{\l}(1+\kappa(r)), \
\text {with} \ |\kappa(r)|\le \frac{1}{10}, \ r=0,\ldots d.
\ee
\el
\pr
The required expression follows from (\ref{eq:der.expression}) with
$$
\kappa=\frac{\beta_1}{\alpha_1t}+\sum^{d-r}_{k=2}\frac{B(k,r)}{\l^{k-1}}\left
(\frac{\alpha_k}{\alpha_1}+\frac{\beta_k}{\alpha_1t}\right ).
$$
Hence, by the assumptions on $t$ and $\lambda$ we have
$$
|\kappa|\le \frac{1}{20}+ \frac{1}{\l}\sum^{d-r}_{k=2}B(k,r)\left (|\frac{\alpha_k}{\alpha_1}|+
\frac{1}{20}\right )\le \frac{1}{10}. \ \square
$$
Thus, we have to count the number of sign changes in sequence (\ref{eq:der.expression1}),
with $r=0,\ldots,d$, for different $\l$.

\smallskip

Assume first that $\alpha_1>0$. In this case there are no sign changes in (\ref{eq:der.expression1}),
and we conclude for each $\l\ge \l_0$ and $t\ge t_0$ we have $\nu_{Q_t}(\l)=0$.
Therefore $Q_t(z)$ does not have real roots on $[\l_1,\infty)$.

\smallskip

Now we consider the case $\alpha_1<0$. Put $c_1(d)=\min_{r=0,\ldots,d-1} \ C(r),
\ c_2(d)=\max_{r=0,\ldots,d-1} \ C(r)$, and define $t_1:=\frac{2\l_0}{c(d)}$, \ $A_1:=\frac{1}{2}c_1(d)|\alpha_1|,
\ A_2:=2c_2(d)|\alpha_1|$.

\smallskip

We get immediately that for each $t\ge t_1$
$$
\frac{C(r)t\alpha_1}{\l_0}(1+\kappa(r))<-1, \ r=0,\ldots d-1,
$$
and for $\l_1(t)=A_1t, \ \l_2(t)=A_2t$
$$
\frac{C(r)t\alpha_1}{\l_1(t)}(1+\kappa(r))<-1, \ \frac{C(r)t\alpha_1}{\l_2(t)}(1+\kappa(r))>-\frac{1}{2}, \ r=0,\ldots d-1.
$$
By Lemma \ref{lem:simple.expr} we get $Q^{(r)}_t(\l_0), Q^{(r)}_t(\l_1(t))<0, \ r=0,\ldots,d-1$, while we have $Q^{(d)}_t\equiv 1$, and hence $\nu_{Q_t}(\l_1(t))=1$. For $\l_2(t)$ we obtain all $Q^{(r)}_t(\l_2(t)), \ r=0,\ldots,d,$ positive, i.e. $\nu_{Q_t}(\l_2(t))=0=\nu_{Q_t}(\infty).$ Therefore there are no real roots of $Q_t(z)$ between $\l_0$ and $\l_1(t)$, there is exactly one real root of $Q_t(z)$ between $\l_1(t)$ and $\l_2(t)$, and no real roots in $[\l_2(t),\infty).$ This completes the proof of Proposition \ref{Prop:real.roots} and of Theorem \ref{thm:nodes.Prony.curves}. $\square$

\smallskip

\noindent{\bf Remark.} The cases $t\to -\infty$ and/or $z\to -\infty$ are reduced
to the case above by the substitutions $\tau=-t$ and $w=-z,$ respectively.

\section{Prony varieties for two nodes}\label{Sec:two.nodes}
\setcounter{equation}{0}

In this section we illustrate some of the results above, providing a complete
description of the Prony varieties in the case of two nodes, i.e. for $d=2$. Put
$\mu=(\mu_0,\mu_1,\mu_2,\mu_3)$.

\smallskip

For $q=0$ the varieties $S_0(\mu)$ are three-dimensional hyperplanes in ${\cal P}_2\cong {\mathbb R}^4$,
defined by the equation $a_1+a_2=\mu_0$.

\smallskip

For $q=1=d-1$ the varieties $S_1(\mu)$ are two-dimensional subvarieties in ${\cal P}_2$, defined by the equations

\be\label{eq:Prony.fol.d.2.1}
a_1+a_2=\mu_0, \ \ a_1x_1+a_2x_2=\mu_1.
\ee
This gives

\be\label{eq:Prony.fol.d.2.3}
a_1=\frac{\mu_0x_2-\mu_1}{x_2-x_1}, \ a_2 = \frac{-\mu_0x_1+\mu_1}{x_2-x_1},
\ee
which is a special case, for $d=2$, of expressions (\ref{eq:ampl.fprmula.yoni}).

\medskip

Consider now the case $q=2=d=2d-2$. Here the varieties $S_2(\mu)$ are (generically) algebraic curves in ${\cal P}_2$,
defined by the equations

\be\label{eq:Prony.fol.d.2.4}
a_1+a_2=\mu_0, \ \ a_1x_1+a_2x_2=\mu_1, \ \ a_1x^2_1+a_2x^2_2=\mu_2,
\ee
For the corresponding curve $S^x_2(\mu)$ in the nodes space
${\cal P}^x_2\cong \Delta_2$ we obtain from Theorem \ref{thm:sing.Sq3}
the equation $\mu_{1}\sigma_1+\mu_{0}\sigma_2 =-\mu_2,$ or
\be\label{eq:Prony.fol.d.2.5}
\mu_0x_1x_2-\mu_1(x_1+x_2)+\mu_2=0.
\ee
This equation leads to three different possibilities:

\smallskip

\noindent 1. If $\mu_0\ne 0$, then the curve $S^x_2(\mu)$ is a hyperbola

\be\label{eq:Prony.fol.d.2.6}
(x_1-\frac{\mu_1}{\mu_0})(x_2-\frac{\mu_1}{\mu_0})+\frac{\mu_0\mu_2-\mu_1^2}{\mu_0^2}=0,
\ee
which is non-singular for $\mu_0\mu_2-\mu_1^2\ne 0$, and degenerates into two orthogonal coordinate lines,
crossing at the diagonal $\{x_1=x_2\},$ for $\mu_0\mu_2-\mu_1^2=0$.

\smallskip

\noindent 2. If $\mu_0=0$, but $\mu_1\ne 0$ then the curve $S^x_2(\mu)$ is a straight line

\be\label{eq:Prony.fol.d.2.7}
x_1+x_2=\mu_2/\mu_1.
\ee
\noindent 3. Finally, if $\mu_0=\mu_1=0$, but $\mu_2\ne 0$ then the curve $S^x_2(\mu)$ is empty,
and for $\mu_0=\mu_1=\mu_2=0$ it coincides with the entire plane ${\cal P}^x_2$.

\medskip

\begin{figure}
	\centering
	\includegraphics[scale=0.6]{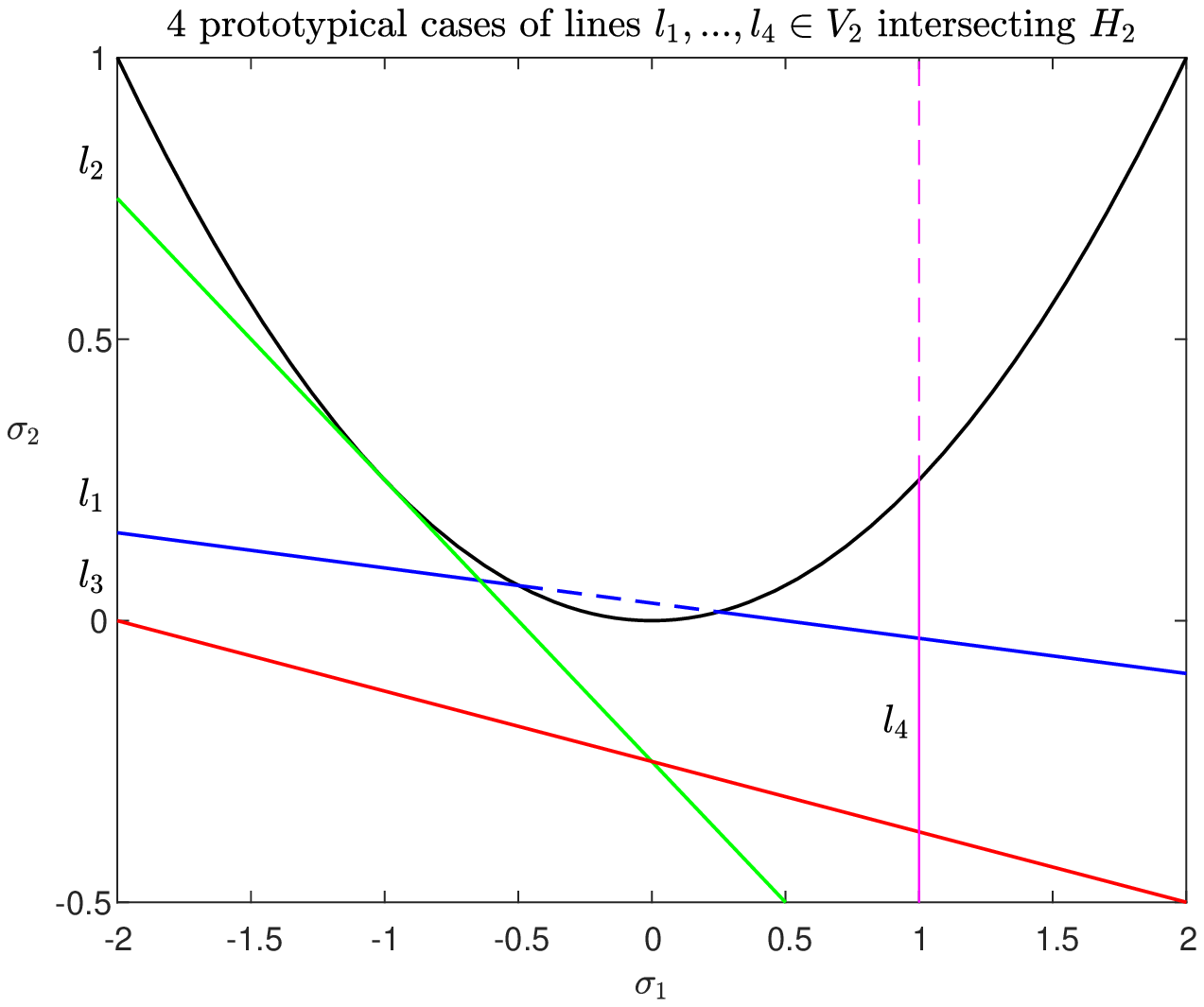}
	\includegraphics[scale=0.75]{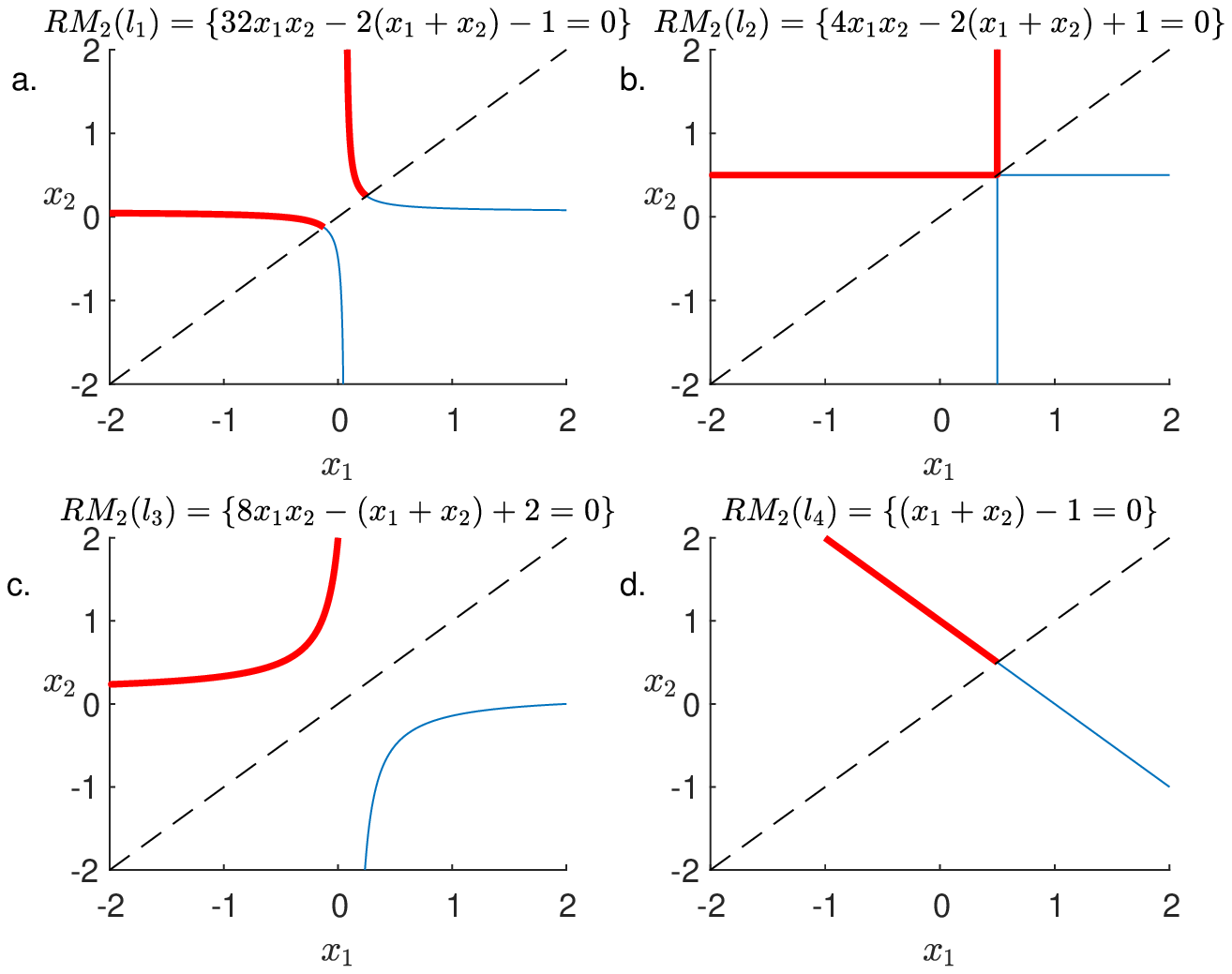}
	\caption{Visualized is the diffeomorphism $RM:H_2\to {\cal P}^x_2$ acting on 4 prototypical lines $l_1,\;l_2\;,l_3\;,l_4 \in {\cal V}_2$
	intersected with $H_2$  (the open set outside the parabola on the upper figure).
	In the bottom figure, the highlighted parts in the subplots a,b,c,d are the images, under $RM_2$, of $l_1,\;l_2\;,l_3\;,l_4$ intersected with $H_2$, respectively.}
	\label{fig.iso}
\end{figure}

It is instructive to interpret the cases (1-3) above in terms of the relative position,
with respect to the set $H_2$ of hyperbolic polynomials $Q$, of the straight line $L_2(\mu)$.
This line is defined in the space ${\cal V}_2$ of polynomials $Q(z)=z^2+\sigma_1 z+\sigma_2$
by system (\ref{eq:Prony.fol.new.4}),
i.e. by the equation $\mu_{1}\sigma_1+\mu_{0}\sigma_2 =-\mu_2$. Figure \ref{fig.iso}
illustrates possible positions of the line $L_2(\mu)$ with respect to the set $H_2$ of hyperbolic polynomials.

\smallskip

The discriminant $\Delta(\sigma_1,\sigma_2)=\sigma^2_1-4\sigma_2$ of $Q(z)=z^2+\sigma_1 z+\sigma_2$
is positive for $Q\in H_2$.
Therefore $H_2$ is the part under the parabola $P=\{\sigma_2=\frac{1}{4}\sigma^2_1\}$ in ${\cal V}_2$.
(Compare Figure \ref{fig.iso}).
The case $\mu_0\ne 0$ corresponds to the lines $L_2(\mu)$, nonparallel to the $\sigma_2$-axis of ${\cal V}_2$.
These lines may cross the parabola $P$ at two points (line $l_1$ on Figure \ref{fig.iso}),
at one point, if tangent to $P$ (line $l_2$ on Figure \ref{fig.iso}),
or they may not cross $P$ at all, and then they are entirely contained in $H_2$ (line $l_3$ on Figure \ref{fig.iso}).
These cases correspond to $\mu_0\mu_2-\mu_1^2<0$, $\mu_0\mu_2-\mu_1^2=0$ and $\mu_0\mu_2-\mu_1^2>0$, respectively.

\smallskip

The corresponding Prony curves $S^x(\mu),$ which are the images under the
root map $RM$ of the lines $L_2(\mu)$ intersected with $H_2$, are shown on the bottom part of Figure \ref{fig.iso}:

\smallskip

For the line $L_2(\mu)$ crossing the parabola $P$ at two points (like $l_1$ on Figure \ref{fig.iso})
the corresponding hyperbola $S^x_2(\mu)$
crosses the diagonal in the plane ${\cal P}^x_2$, i.e. it contain a collision of the nodes $x_1,x_2$.

Notice that the Prony curve $S^x(\mu)$ remains non-singular at the crossing point. This fact holds also for 
the general case of ``double collisions'' on the Prony curves, and we plan to present it separately.

\smallskip

For the line $L_2(\mu)$ tangent to the parabola $P$ (like $l_2$ on Figure \ref{fig.iso}),
the corresponding hyperbola $S^x_2(\mu)$
degenerates into two orthogonal coordinate lines, crossing at a certain point on the diagonal $\{x_1=x_2\},$.

\smallskip

For the line $L_2(\mu)$ entirely contained in $H_2$ (like $l_3$ on Figure \ref{fig.iso})
the corresponding hyperbola $S^x_2(\mu)$ does not cross the diagonal $\{x_1=x_2\},$
and so it does not lead to the nodes collision.

\medskip

For $\mu_0=0$, but $\mu_1\ne 0$, the lines $L_2(\mu)$ are parallel to
the $\sigma_2$-axis of $V_2$ (like $l_4$ on Figure \ref{fig.iso}).
They cross the parabola $P$ at exactly one point. The corresponding
curve $S^x_2(\mu)$ for $l_4$ is a straight line $x_1+x_2=-\frac{\mu_2}{\mu_1}$.

\smallskip

\bibliographystyle{myplain}
\bibliography{bib}{}

\end{document}